\documentclass[leqno,12pt]{article}
\usepackage{latexsym}
\usepackage{amssymb}
\usepackage{amsmath}
\usepackage{amssymb}

\usepackage{graphicx}
\usepackage{epsfig}

 0

\newcommand{\erw}[1]{\mbox{E} [#1] }
\newcommand{\ta}{\tilde{a}_n}
\newcommand{\tb}{\tilde{b}_n}

\renewcommand{\theequation}{\thesection.\arabic{equation}}
\makeatletter\@addtoreset{equation}{section}\makeatother

\newcommand{\ba}{\begin{array}}
\newcommand{\ea}{\end{array}}
\newcommand{\beqohne}{\begin{eqnarray*}}
\newcommand{\eeqohne}{\end{eqnarray*}}
\newcommand{\beohne}{\begin{equation*}}
\newcommand{\eeohne}{\end{equation*}}

\def\3{\ss}

\newcommand{\bea}{\begin{eqnarray*}}
\newcommand{\eea}{\end{eqnarray*}}
\newcommand{\be}{\begin{eqnarray}}
\newcommand{\ee}{\end{eqnarray}}
\newcommand{\beq}{\begin{equation}}
\newcommand{\eeq}{\end{equation}}

\begin{document}

\title{Some asymptotic properties of the spectrum of the Jacobi ensemble }

\author{
{\small Holger Dette} \\
{\small Ruhr-Universit\"at Bochum} \\
{\small Fakult\"at f\"ur Mathematik} \\
{\small 44780 Bochum, Germany} \\
{\small e-mail: holger.dette@rub.de}\\
{\small FAX: +49 234 3214 559}\\
\and
{\small Jan Nagel} \\
{\small Ruhr-Universit\"at Bochum} \\
{\small Fakult\"at f\"ur Mathematik} \\
{\small 44780 Bochum, Germany} \\
{\small e-mail: jan.nagel@rub.de}\\
}

\maketitle

\begin{abstract}
\noindent For the random eigenvalues with density corresponding to the Jacobi ensemble $$c \cdot \prod_{i < j} \ | \ \lambda_i - \lambda_j \ |^\beta \prod^n_{i=1} (2 - \lambda_i)^a (2 + \lambda_i)^b I_{(-2,2)} (\lambda_i ) $$ $(a, b > -1, \beta > 0) $
a strong uniform approximation by the roots of
the Jacobi polynomials is derived if the parameters $a, b,$ $\beta$ depend on $n$ and $n \to \infty$. Roughly speaking, the eigenvalues can be
uniformly approximated by roots of Jacobi polynomials with parameters $((2a+2)/\beta -1, (2b+2)/\beta-1)$, where the error is of order $\{ \log
n/(a+b) \}^{1/4}$. These results are used to investigate the asymptotic properties of the corresponding spectral distribution if $n \to \infty$
and the parameters $a, b$ and $\beta$ vary with $n$. We also discuss further applications in the context of multivariate random $F$-matrices.

\end{abstract}



\medskip

\section{Introduction}
\def\theequation{1.\arabic{equation}}
\setcounter{equation}{0}

The three classical ensembles of random matrix theory are the Hermite, Laguerre and Jacobi ensembles. The Hermite or Gaussian ensembles arise
in physics and are obtained as the eigenvalue distribution of a symmetric matrix with Gaussian entries. The Laguerre or Wishart ensembles
appear in statistics and correspond to the distribution of the singular values of a Gaussian matrix. Similarly, the Jacobi ensembles are
objects of statistical interest and are motivated by multivariate analysis of variance (MANOVA; see Muirhead (1982)). Associated with each
ensemble there is a real positive parameter $\beta$ which is usually considered for three values. The case $\beta = 1$ corresponds to real
matrices, while the ensembles for $\beta = 2$ and $\beta = 4$ arise from complex and quaternion random matrices, respectively, according to
Dyson's (1962) threefold classification. Dyson also observed that the eigenvalue distributions correspond to the Gibbs distribution for the
classical Coloumb gas at three different temperatures. In other words -- starting from the physical interpretation -- for many decades it was
only known that there exist random matrix models for the Coloumb gas at three  different temperatures. Recently Dumitriu and Edelman (2002)
 provided tridiagonal random matrix models for the $\beta$-Hermite and $\beta$-Laguerre ensembles
for all $\beta > 0$. Dette and Imhof (2007) derived strong uniform approximations of the random eigenvalues of these ensembles by  the
(deterministic) roots of Hermite and Laguerre polynomials. The development of a matrix model corresponding to the Jacobi ensemble was an open
problem, which was recently considered by Lippert  (2003) and  Killip and Nenciu (2004) and Edelman and Sutton (2006). In particular, Killip
and Nenciu (2004) found a tridiagonal matrix model for the $\beta$-Jacobi ensemble with density
\be \label{1.1} f(\lambda ) = c \prod_{i < j} | \lambda _i - \lambda _j |^{\beta } \prod_{i=1}^{n} (2-\lambda _i)^a(2+\lambda _i)^b  I_{(-2,2)} (\lambda_i ) , \ee $(a, b > -1, \beta > 0)$, where
the entries are simple functions of independent random variables with a beta distribution on the interval $[-1,1]$ [see equations (\ref{2.2})
and (\ref{2.3}) in Section 2 for more details]. It is the purpose of the present paper to obtain further insight in the stochastic properties
of the random eigenvalues with density given by (\ref{1.1}). In Section 2 we introduce the triangular matrix proposed by Killip and Nenciu
(2004) and present a uniform approximation of the random eigenvalues with density (\ref{1.1}) by roots of Jacobi polynomials. Roughly speaking,
if $n \to \infty$, the random eigenvalues can be uniformly approximated by roots of the Jacobi polynomials $P_n^{((2a+2)/\beta-1,
(2b+2)/\beta-1)} (x/2)$, where the error of this approximation is of order $\{ \log n/(a+b) \}^{1/4}$. These results are used in Section 3 to
investigate the asymptotic properties of the spectrum if the parameters $a$, $b$ and $\beta$ vary simultaneously with $n$ and $n \to \infty$.
In Section 4 we study as an application the eigenvalue distribution of a multivariate $F$-matrix which was also investigated by Silverstein
(1985b) and Collins (2005). We present several extensions of these results. In particular, we consider almost sure convergence [Silverstein
(1985b) discussed convergence in probability while Collins (2005) considered the expectation of the empirical spectral distribution] and
investigate the case, where the parameters and the sample size converge to infinity with different order [Silverstein (1985b) and Collins
(2005) discussed the case where $a \sim \gamma_1 n, \  b \sim \gamma_2 n$]. Finally, some technical results are given in an appendix.

\section{Strong approximation  of the Jacobi ensemble}
\def\theequation{2.\arabic{equation}}
\setcounter{equation}{0}

Recall the definition of the Jacobi ensemble in (\ref{1.1}) and let for $p,q>0$ $\alpha \sim B (p,q)$ denote a Beta-distributed random variable
on the interval $[-1, 1]$ with density
\be \label{2.1}  2^{1-p-q} \frac{\Gamma (p+q)}{\Gamma (p)\Gamma (q)} (1-x)^{p-1}(1+x)^{q-1}
I_{(-1,1)} (x). \ee
It was shown by Killip and Nenciu (2004) that if $\alpha_0, \alpha_1, \dots, \alpha_{2n-2}$ are independent random
variables
with distribution \be \label{2.2} \alpha _k \sim
\begin{cases}
B(\tfrac{2n-k-2}{4} \beta +a+1, \tfrac{2n-k-2}{4} \beta +b+1) & \mbox{ if }
k \mbox{ even,} \\
B(\tfrac{2n-k-3}{4} \beta +a+b+2, \tfrac{2n-k-1}{4} \beta ) & \mbox{ if } k
\mbox{ odd,}
\end{cases}
\ee then the joint density of the (real) eigenvalues of the tridiagonal matrix \be \label{2.3} J := \begin{pmatrix}
                b_1 & a_1    &         &         \\
                a_1 & b_2    & \ddots  &         \\
                    & \ddots & \ddots  & a_{n-1} \\
                    &        & a_{n-1} & b_n
        \end{pmatrix}
\ee
with entries
\begin{eqnarray*}
b_{k+1} &=& (1-\alpha _{2k-1} )\alpha _{2k} - (1+\alpha _{2k-1} )\alpha _{2k-2} \ ,\\
a_{k+1} &=& \left \{ (1- \alpha _{2k-1} )(1-\alpha _{2k} ^2)(1+\alpha _{2k+1} ) \right \} ^{1/2}
\end{eqnarray*}
$(\alpha_{2n-1} = \alpha_{-1} = \alpha_{-2} = -1)$  is given by the Jacobi ensemble (\ref{1.1}).
In the following discussion we consider the asymptotic properties of the eigenvalues of the random matrix (\ref{2.3})
[or equivalently of the Jacobi ensemble (\ref{1.1})],
where $n \to \infty$ and the parameters $a$, $b$ and $\beta$ in (\ref{1.1}) also vary with  $n$. An important tool of
our analysis are the Jacobi polynomials $P_n^{(\gamma, \delta)}(x)$, which are defined as the polynomials of degree $n$
with leading coefficient $(n + \gamma + \delta + 1)_n / 2^n n!$ satisfying the orthogonality relation
\be \label{2.5}
 \int_{-1}^1 P_n^{(\gamma,\delta)}(x) P_m^{(\gamma,\delta)}(x)
(1-x)^\gamma(1+x)^\delta dx = 0 \qquad \mbox{if} \quad m\neq n \ee [see Szeg\"o (1975)].
Here and throughout this paper $(a)_n = \Gamma (a+n)/\Gamma (a) $ denotes the Pochhammer symbol.
The following results provide an almost sure uniform approximation
of the random eigenvalues of the Jacobi ensemble by the roots of  orthogonal polynomials, if $n \to \infty$. We begin with a statement of an
exponential bound for the probability of a maximal deviation between the roots of the Jacobi polynomials $P_n^{((2a+2)/\beta-1, (2b+2)/\beta-1)}
(x/2)$ and random eigenvalues of the Jacobi ensemble.

\bigskip

{\bf Theorem 2.1.} {\it Let $\lambda_1^{(n)} \leq \dots \leq \lambda_n^{(n)}$ denote the ordered eigenvalues
with density given by the Jacobi ensemble
(\ref{1.1}) and $x_1^{(n)} < \dots < x_n^{(n)}$ denote the ordered roots of the Jacobi polynomial
$$
P_n^{(\frac{2a+2}{\beta}-1, \frac{2b+2}{\beta}-1)} (\tfrac {1}{2}x) \: ,
$$
then the following inequality holds for any  $\varepsilon \in (0,1]$
\be  \nonumber
 P\left( \underset{1 \leq j \leq n}{\max } |\lambda _j^{(n)} - x_j^{(n)} | > \varepsilon \right)
                    \leq  4  (2n-1) \exp \left \{ \left( \log ( 1+\frac{\varepsilon ^2}{648+2\varepsilon ^2} )
                     - \frac{\varepsilon ^2}{648+2\varepsilon ^2} \right) (a+b+2) \right \}\: .  \\
\ee }

\bigskip

{\bf Proof of Theorem 2.1.} Interchanging the rows and columns of the matrix $J$ defined in (\ref{2.3}) it follows that the matrix
\be \label{2.7}
 \tilde{J} := \begin{pmatrix}
                {b'}_1 & {a'}_1    &         &         \\
                {a'}_1 & {b'}_2    & \ddots  &         \\
                    & \ddots & \ddots  & {a'}_{n-1} \\
                    &        & {a'}_{n-1} & {b'}_n
        \end{pmatrix}
\ee with entries
\begin{eqnarray} \label{2.8}
{b'}_{k+1} &:=& b_{n-k} =  (1-\alpha _{2n-2k-3} )\alpha _{2n-2k-2} -
(1+\alpha _{2n-2k-3} )\alpha _{2n-2k-4} \ , \\
{a'}_{k+1} &:=& a_{n-k-1} = \left \{ (1- \alpha _{2n-2k-5} )(1-\alpha
_{2n-2k-4} ^2)(1+\alpha _{2n-2k-3} ) \right \} ^{1/2}
\label{2.9}
\end{eqnarray}
has the same eigenvalues as the matrix $J$. This implies that the joint density of the eigenvalues of the matrix $\tilde J$ is also given by
the $\beta$-Jacobi ensemble defined in (\ref{1.1}). We now consider the deterministic matrix \be \label{2.10}
 D := \begin{pmatrix}
                d_1 & c_1    &         &         \\
                c_1 & d_2    & \ddots  &         \\
                    & \ddots & \ddots  & c_{n-1} \\
                    &        & c_{n-1} & d_n
        \end{pmatrix} \: ,
\ee where we essentially replace the random variables in (\ref{2.8}) and (\ref{2.9}) by their corresponding expectations, that is
\begin{eqnarray*}
d_{k+1} & := &  (1-\erw{\alpha _{2n-2k-3}}
)\erw{\alpha _{2n-2k-2}} - (1+\erw{\alpha _{2n-2k-3}} )\erw{\alpha
_{2n-2k-4}} \ , \\
c_{k+1} & := & \left \{ (1- \erw{\alpha _{2n-2k-5}} )(1-\erw{\alpha
_{2n-2k-4}} ^2 )(1+\erw{\alpha _{2n-2k-3}} ) \right \} ^{1/2} \ .
\end{eqnarray*}
A straightforward calculation [observing that the expectation 
of a random variable with density (\ref{2.1}) is given by
$(q-p)/(p+q)$
]yields
\begin{eqnarray}
d_{k+1} &=&  \frac{2(\tilde b ^2 - \tilde a ^2)}{(2k+\tilde a +\tilde b )(2k+\tilde a +\tilde b +2)} \: ,  \qquad k = 0, \dots, n-2 \: ;
\label{2.11a}
\\
 c_{k+1} &=& \frac{4}{2k +\tilde a + \tilde b +2} \left\{ \frac{(k+\tilde a
+\tilde b +1)(k + \tilde a +1)
(k+\tilde b +1)(k+1)}{(2k + \tilde a +\tilde b +3)(2k +\tilde a +\tilde b
+1)} \right\} ^{1/2}\label{2.12}
\end{eqnarray}
$ (k = 0, \dots, n-3)$
where $\tilde a = \frac{2a+2}{\beta} $, $\tilde b = \frac{2b+2}{\beta} $,
and (observing that $\alpha_{-1} = \alpha_{-2} = -1$)
\begin{eqnarray}
d_n &=& \frac{2(\tilde b -\tilde a )}{2n +\tilde a +\tilde b -2} \label{2.13} \ ,\\
c_{n-1} &=&  \frac{4}{2n+\tilde a +\tilde b -2} \left\{ \frac{(n+\tilde a -1)(n+\tilde b -1)(n-1)}{(2n +\tilde a +\tilde b -3)} \right\} ^{1/2}
\: .\label{2.14}
\end{eqnarray}
For the calculation of the eigenvalues of the matrix $D$ we consider the determinant det$(x \ I_n - D)$, then it follows by an expansion  with respect to the last
row that \be \label{2.15} \det (x \ I_n -D) = (x-d_n) \ G_{n-1} (x) - c^2_{n-1} \ G_{n-2} (x) \: , \ee
where the polynomials $G_0 (x), \dots, G_{n-1} (x)$ are
defined recursively by the three term recurrence relation \be \label{2.16} G_k (x) = (x- d_k) \ G_{k-1} (x) - c^2_{k-1} \ G_{k-2} (x) \ee 
$(k=1,\dots,n-1 ; \ \  G_0 (x) :=1, \ \ G_{-1}(x) := 0)$. Now a straightforward calculation and  a comparison with the three term recurrence
relation for the monic Jacobi polynomials [see e.g.\ Chihara (1978), p.\ 220] yields \be \label{2.17} G_k (x) = 2^k \hat P_k^{(\tilde a,\tilde b)}
(\tfrac{1}{2} {x}) ; \qquad k = 0, \dots, n-1, \ee where $\hat P_k^{(\tilde a, \tilde b)} (x)$ denotes the $k$th monic Jacobi polynomial, i.e. \be \label{2.18} \hat
P_k^{(\tilde a,\tilde b)} (x) = \frac {2^k \ k!}{(k+\tilde a+\tilde b+1)_k} \ P^{(\tilde a,\tilde b)}_k (x) \: . \ee Combining equations (\ref{2.16}), (\ref{2.17}) and (\ref{2.18}) we
obtain by a tedious but straightforward calculation
\begin{eqnarray*}
\det (I_n \ x - D) &=& G_n (x) + (s_n - d_n) \ G_{n-1} (x) + (r^2_{n-1} - c^2_{n-1}) \ G_{n-2} (x) \\
&& \\
&=& \frac{4^n n!}{(n+\tilde a+\tilde b+1)_n} \Bigl\{ P_{n}^{(\tilde a ,\tilde b
)}(\tfrac{1}{2} x) -
\frac{(\tilde b -\tilde a )(2n +\tilde a +\tilde b -1)}{(2n +\tilde a
+\tilde b -2)(n +\tilde a +\tilde b )}
            P_{n-1}^{(\tilde a ,\tilde b )}(\tfrac{1}{2} x)  \\
    & &  - \frac{(2n +\tilde a +\tilde b )(n+\tilde a -1)(n+\tilde
b -1)}{(2n+\tilde a +\tilde b -2)(n+\tilde a +\tilde b )(n+\tilde a +\tilde
b -1)} P_{n-2}^{(\tilde a ,\tilde b )}(\tfrac{1}{2} x) \Bigr\}  \\
&& \\
    &=& \frac{4^n n!}{(n+\tilde a+\tilde b+1)_n} \Bigl\{ P_{n}^{(\tilde a ,\tilde b
)}(\tfrac{1}{2} x) -
\frac{n+\tilde b }{n+\tilde a +\tilde b } P_{n-1}^{(\tilde a ,\tilde b
)}(\tfrac{1}{2} x)  \\
        && ~~~~~~~~
+ \frac{(2n +\tilde a +\tilde b )(n+\tilde a -1)}{(n + \tilde a +\tilde b )(n+\tilde a +\tilde b -1)} P_{n-1}^{(\tilde a -1,\tilde b
)}(\tfrac{1}{2} x) \Bigr\} \: ,
\end{eqnarray*}
where we have used the identity \be \label{2.19} (n+\tilde b -1) P_{n-2}^{(\tilde a ,\tilde b )}(x) = (n +\tilde a +\tilde b -1)
P_{n-1}^{(\tilde a ,\tilde b )}(x) - (2n +\tilde a +\tilde b -2) P_{n-1}^{(\tilde a -1,\tilde b )}(x) \ee in the last step [see Abramovich and
Stegun (1965), equation (22.7.18)]. A further application of this identity to the second polynomial yields
\begin{eqnarray*}
    \mbox{det}(I_nx-D)
&=& \frac{4^n n!}{(n+\tilde a+\tilde b+1)_n} \Bigl\{\frac{2n+\tilde a +\tilde b
}{n+\tilde a +\tilde b } P_{n}^{(\tilde a -1,\tilde b )}(\tfrac{1}{2} x)  \\
&& \\
    & & \qquad \
 + \frac{(2n +\tilde a +\tilde b )(n+\tilde a -1)}{(n + \tilde a +\tilde b
)(n+\tilde a +\tilde b -1)} P_{n-1}^{(\tilde a -1,\tilde b )}(\tfrac{1}{2}
x) \Bigr\} \\
&& \\
    &=& \frac{4^n n!}{(n+\tilde a +\tilde b +1)_n} \Bigl\{
\frac{2n+\tilde a +\tilde b }{n+\tilde a +\tilde b } P_{n}^{(\tilde a
-1,\tilde b )}(\tfrac{1}{2} x)
 - \frac{2n+\tilde a +\tilde b }{n+\tilde a +\tilde b } P_{n}^{(\tilde a
-1,\tilde b )}(\tfrac{1}{2} x)          \\
    & & \quad + \frac{(2n +\tilde a +\tilde b )(2n +\tilde a +\tilde b
-1)}{(n + \tilde a +\tilde b )(n + \tilde a +\tilde b -1)}
P_{n}^{(\tilde a -1,\tilde b -1)}(\tfrac{1}{2} x) \Bigr\}   \\
&& \\
    &=& \frac{4^n n!}{(n+(\tilde a -1)+(\tilde b -1)+1)_n}
P_{n}^{(\tilde a -1,\tilde b -1)}(\tfrac{1}{2} x) \: ,
\end{eqnarray*}
where we have used the identity \be \label{2.20} ~~~~~~~~ (n +\tilde a -1) P_{n-1}^{(\tilde a -1 ,\tilde b )}(x) = (2n+\tilde a + \tilde b
-1) P_{n}^{(\tilde a -1,\tilde b -1)}(x) - (n+\tilde a +\tilde b -1) P_{n}^{(\tilde a -1,\tilde b )}(x) \ee for the second equality
[see Abramovich and
Stegun (1965), equation (22.7.19)].
Consequently, the eigenvalues of the matrix $D$ are given by the roots $x_1^{(n)} < \dots < x_n^{(n)}$ of the Jacobi polynomial $P_n^{(\tilde a - 1, \
\tilde b - 1)} (\frac {1}{2} x)$. A similar argument as in Silverstein (1985a) now shows that \be \label{2.21} \underset{1 \leq j \leq n}{\max
} |\lambda _j^{(n)} - x_j^{(n)} | \leq \underset{1 \leq k \leq n}{\max } \sum_{j=1}^n | \tilde{J} _{kj} -D_{kj} |
                \leq 4 \left \{ 3X_n \right \} ^{1/2} +
6X_n  \: , \ee where the elements of the matrices $\tilde J$ and $D$ are denoted by $\tilde J_{ij}$ and $D_{ij}$, respectively, and the random
variable $X_n$ is defined by \be \label{2.22} X_n := \underset{0 \leq k \leq 2n-2}{\max } |\alpha _{k} - \erw{\alpha _{k} } |  \: . \ee This implies for the
probability in Theorem 2.1
\begin{eqnarray}\label{2.23}
P\left( \underset{1 \leq j \leq n}{\max } |\lambda _j^{(n)} - x_j^{(n)} | > \varepsilon
\right)
&\leq & P\left( 2\{ 3X_n \} ^{1/2} + 3X_n > \tfrac{\varepsilon }{2} \right)
     \leq
P\left( 3\{ 3X_n \} ^{1/2} > \tfrac{\varepsilon }{2} \right)  \\
&=&  P\left( X_n > \tfrac{\varepsilon ^2}{108} \right)
    \leq \sum_{k=0}^{2n-2} P\left( |\alpha _k - \erw{\alpha _k} | >
\tfrac{\varepsilon ^2}{108} \right) \nonumber
\end{eqnarray}
whenever $\varepsilon \in (0,1]$. Observing Lemma A.1 in the Appendix and that $\alpha_k = 1 - 2 \beta_k$, where
$\beta_k$ is the corresponding Beta distribution on the interval $[0,1]$, it follows that
\begin{eqnarray*}
P \ \Bigl ( \max_{1 \leq j \leq n} \mid \lambda_j^{(n)} - x_j^{(n)} \mid \ > \varepsilon \Bigr ) & \leq & 4 \sum^{2n-2}_{k=0} \exp
\left ( c \ ( \frac {2n-2-k}{2} \ \beta + a+b+2) \right ) \\
& \leq & 4 (2n-1) \exp  \left ( c \ (a+b+2) \right ) \
\end{eqnarray*}
where the constant $c = c(\varepsilon)$ is given by \be \label{2.24} c = \log \left ( 1 + \frac {\varepsilon^2}{648 + 2 \varepsilon^2} \right )
- \frac {\varepsilon^2}{648 + 2 \varepsilon^2}. \ee This proves the assertion of the theorem. \hfill $\Box$

\bigskip

Note that the constant $c$ in (\ref{2.24}) is negative, and Theorem 2.1 therefore indicates that the random eigenvalues
of the Jacobi ensemble can be approximated by the deterministic roots of the Jacobi polynomial $P_n^{(\tilde a-1,\ \tilde b-1)} (\frac {1}{2}x)$
with a high probability if $a+b$ is large. The following result makes this statement more precise and provides a strong uniform
approximation of the eigenvalues of the Jacobi ensemble by the roots of the Jacobi polynomial $P_n^{(\tilde a-1,\ \tilde b-1)} (\frac {1}{2}x)$, if the parameters in
(\ref{1.1}) converge sufficiently fast to infinity. The proof follows by similar arguments as the proof of Theorem 2.2 in Dette and
Imhof (2007) and is therefore omitted.

\bigskip

{\bf Theorem 2.2.} {\it Let $\lambda_1^{(n)} \leq \dots \leq \lambda_n^{(n)}$ denote the ordered random eigenvalues with density given by the Jacobi ensemble (\ref{1.1})
with parameters $a = a_n, b = b_n, \beta = \beta_n$ and $x_1^{(n)} < \dots < x_n^{(n)}$ denote the ordered roots of the Jacobi polynomial
$P_n^{(\tilde a_n - 1, \ \tilde b_n - 1)} (\frac {1}{2}x)$ where \be \label{2.25} \ta = \tfrac{2a_n+2}{\beta_n} , \qquad \tb = \tfrac{2b_n+2}{\beta_n}.
\ee If \be \label{2.26} \lim _{n \rightarrow \infty } \frac{a_n+b_n}{\log n} = \infty \: , \ee then the inequality
\[
\max_{1\leq j\leq n}\left|\lambda_j^{(n)}-x_j^{(n)}\right| \leq \left(\frac {\log n}{a_n + b_n}\right)^{1/4}S
\]
holds for all $n \geq 2$, where $S$ denotes an a.s.\ finite random variable. In particular, if
\begin{equation}\label{ass4}
\liminf_{n\to\infty}\frac{a_n+b_n}n > 0,
\end{equation}
then there exists an  a.s.\ finite random variable  $S'$ such that the inequality
\[
\max_{1\leq j\leq n}\left|\lambda_j^{(n)}-x_j^{(n)}\right|
\leq
\left(\frac {\log n}n\right)^{1/4}S'
\] holds a.s.\ for all $n\geq 2$. }

\section{Asymptotic spectral properties}
\def\theequation{3.\arabic{equation}}
\setcounter{equation}{0}

In this section we will apply Theorem 2.2 to derive the asymptotic properties of the empirical spectral distribution  \be \label{3.1} F_n^J(\xi ) :=
\frac{1}{n} \sum_{i=1}^n I\{ \lambda _i^{(n)} \leq \xi \} \: , \ee where $\lambda_1^{(n)} \leq \dots \leq \lambda_n^{(n)}$ denote the ordered
eigenvalues of the Jacobi ensemble defined by (\ref{1.1}) with parameters $a = a_n, b = b_n$ and $\beta = \beta_n$.
The results of Section 2 indicate that the
empirical distribution function in (\ref{3.1}) should exhibit a similar asymptotic behaviour as the empirical distribution function \be \label{3.2}
F_n^P(\xi ) := \frac{1}{n} \sum_{i=1}^n I\{ x_i^{(n)} \leq \xi \} \ee of the ordered roots $x^{(n)}_1 < \dots < x^{(n)}_n$ of the Jacobi polynomial $P_n^{(\tilde a_n - 1, \ \tilde b_n
- 1)} (\frac {1}{2}x)$, where $\tilde a_n = (2a_n +2)/\beta_n ,\ \tilde b_n = (2b_n +2)/\beta_n$.
The asymptotic zero distribution of Jacobi polynomials has been studied by several authors [see e.g.\ Gawronski and Shawyer (1991),
Elbert, Laforgia and Rodono (1994),  Dette and Studden (1995) or Kuijlaars and Van Assche (1999) among many others],
and we can use these results and Theorem 2.2 to derive the asymptotic properties of the spectrum of the Jacobi ensemble.  The
following result makes this statement more precise.

\bigskip

{\bf Theorem 3.1.} {\it Consider the empirical distribution functions of the eigenvalues of the Jacobi ensemble (\ref{1.1}) and
the roots of the Jacobi polynomial $P_n^{(\tilde a_n-1, \tilde b_n-1)} (x/2)$ defined by (\ref{3.1}) and (\ref{3.2}), respectively, and let $(\delta_n)_{n \in
\mathbb{N}}, \ (\varepsilon_n)_{n \in \mathbb{N}}$ denote real sequences with $\delta_n > 0$ such that the limit
$$
\lim _{n \rightarrow \infty } F_n^P(\delta _n \xi + \varepsilon _n) = F(\xi )
$$
exists at every coninuity point $\xi$ of $F$. If the conditions
\be \label{3.3}
 \frac{a_n+b_n}{\log n} \xrightarrow[n \rightarrow \infty ]{} \infty \ ,
\qquad \frac{\delta _n^4(a_n+b_n)}{\log n} \xrightarrow[n \rightarrow
\infty ]{} \infty
\ee
are satisfied, then
$$
\lim _{n \rightarrow \infty } F_n^J(\delta _n \xi + \varepsilon _n) = F(\xi )
$$
almost surely at every continuity point $\xi$ of $F$.}

\bigskip

{\bf Proof.} Let $G^P_n$ and $G^J_n$ denote the empirical distribution functions of the rescaled roots
$$
\frac {x_1^{(n)} - \varepsilon_n}{\delta_n}, \dots, \frac {x_n^{(n)} - \varepsilon_n}{\delta_n}
$$ corresponding to the Jacobi polynomial $P_n^{(\tilde a_n - 1, \tilde b_n -1)} (\frac {1}{2}x)$ and of the
eigenvalues $$\frac {\lambda_1^{(n)} - \varepsilon_n}{\delta_n}, \dots, \frac {\lambda_n^{(n)} - \varepsilon_n}{\delta_n}$$ corresponding to the
Jacobi ensemble defined by (\ref{1.1}), respectively. The arguments presented in the proof of Theorem 2.1 show that $G^P_n$ and $G^J_n$ are the empirical distribution
functions of the eigenvalues of the matrices
\begin{eqnarray*}
A_n &:= &
\frac{1}{\delta_n} (D-\epsilon_n I_n),\\
B_n &:= &
\frac{1}{\delta_n} (J-\epsilon_n I_n),
\end{eqnarray*}
respectively. Observing Lemma 2.3 in Bai (1999) we obtain for the Levy-distance between the distribution functions $G^P_n$ and
$G^J_n$
$$
L^3(G_n^P,G_n^J) \leq   \delta_n^{-2} \frac{1}{n} \sum_{i=1}^n \mid \lambda_i^{(n)} -x_i^{(n)}
\mid ^2 \\
\leq \delta_n^{-2} (\underset{1 \leq i \leq n}{\max } |\lambda _i^{(n)} - x^{(n)}_i | )^2\leq S \Bigl ( \Bigl( \frac {\log n}{(a_n +
b_n) \delta^4_n} \Bigr)^{1/2} \Bigr) \: ,
$$
where $S$ denotes an a.s.\ finite random variable. Consequently, we obtain from the assumptions
$$
L(G_n^P,G_n^J) \xrightarrow[n \rightarrow \infty ]{a.s.} 0 \ ,
$$
and the assertion  of Theorem 3.1 follows observing the identities
\begin{eqnarray*}
G^P_n (\xi) &=& F^P_n (\delta_n \xi + \varepsilon_n)\: , \\
G^J_n (\xi) &=& F^J_n (\delta_n \xi + \varepsilon_n) \: .
\end{eqnarray*}
\hfill $\Box$

\bigskip

{\bf Theorem 3.2.} {\it Let $\lambda^{(n)}_1 \leq \dots \leq \lambda^{(n)}_n$ denote the ordered random eigenvalues with density
given by the Jacobi ensemble (\ref{1.1}). If the assumptions of Theorem 3.1 are satisfied and that there exist
 constants $a_1, a_2 \in \mathbb{R} \
\ \
b_1, b_2 \in \mathbb{R}^+$ such that
\begin{align*}
&(i)&    \lim_{n \rightarrow \infty} \frac{1}{\delta_n} \left( \frac{n+\tb
-1}{2n+\ta +\tb -2} -\varepsilon_n \right) &= \frac{a_1}{2} \: ;  \\
&(ii)&   \lim_{n \rightarrow \infty} \frac{1}{\delta_n} \left(
        \frac{n(n+ \ta -1) +(n+\tb -1)(n+\ta +\tb -2)}{(2n+\ta +\tb
-2)^2} -\varepsilon_n \right) &= \frac{a_2}{2} \: ; \\
&(iii)& \lim_{n \rightarrow \infty} \frac{1}{\delta_n^2}
\frac{(n+\tb -1)(n+\ta-1 )n}{(2n+\ta + \tb -2)^3} &= \frac{b_1}{4} \: ;
\\
&(iv) & \lim_{n \rightarrow \infty} \frac{1}{\delta_n^2}
\frac{(n+\tb -1)(n+\ta -1)(n+\ta +\tb -2)n}{(2n+\ta +\tb -2)^4}
                        &= \frac{b_2}{4} \: ;
\end{align*}
then the empirical distribution of the scaled eigenvalues
$$
\frac {\lambda_1^{(n)} - 2(2 \varepsilon_n-1)}{2 \delta_n}  \:, \ \dots, \ \frac {\lambda_n^{(n)} - 2 (2\varepsilon_n-1)}{2 \delta_n}
$$
converges almost surely to a non-degenerate  distribution function, i.e.
$$
\lim_{n \rightarrow \infty} F_n^J(2\delta_n \xi +2(2\varepsilon_n-1))
\stackrel{a.s.}{=} \int_{a_2-2\sqrt{b_2}}^{\xi }
f^{(a_1,a_2,b_1,b_2)}(x)dx \ ,
$$
where
\begin{align*}
f^{(a_1,a_2,b_1,b_2)}(x) =
    \begin{cases}
    \frac{b_1}{2\pi}
\frac{\sqrt{4b_2-(x-a_2)^2}}{(b_2-b_1)x^2+(b_1a_2+b_1a_1-2b_2a_1)
x +b_2{a}_1^2-a_1a_2b_1+{b}_1^2}
&\mbox{if } \mid x-a_2 \mid \leq 2\sqrt{b_2} \ , \\
    0
&\mbox{else.}
    \end{cases}
\end{align*}
}

\bigskip

{\bf Proof.} By Theorem 2.2 in Dette and Studden (1995) it follows that the empirical distribution of the roots of the
Jacobi polynomial $P_n^{(\tilde a_n - 1, \ \tilde b_n - 1)} (\frac {1}{2}x)$ has a non-degenerate limit, that is
$$
 F_n^P(2\delta_n \xi +2(2\varepsilon_n-1))
                    \xrightarrow[n \rightarrow \infty
]{} \int_{a_2-2\sqrt{b_2}}^{\xi } f^{(a_1,a_2,b_1,b_2)}(x)dx \ .
$$
The assertion is now an immediate consequence of Theorem 3.1 \hfill $\Box$

\bigskip

\begin{figure}
\includegraphics[width=8cm,height=7cm]{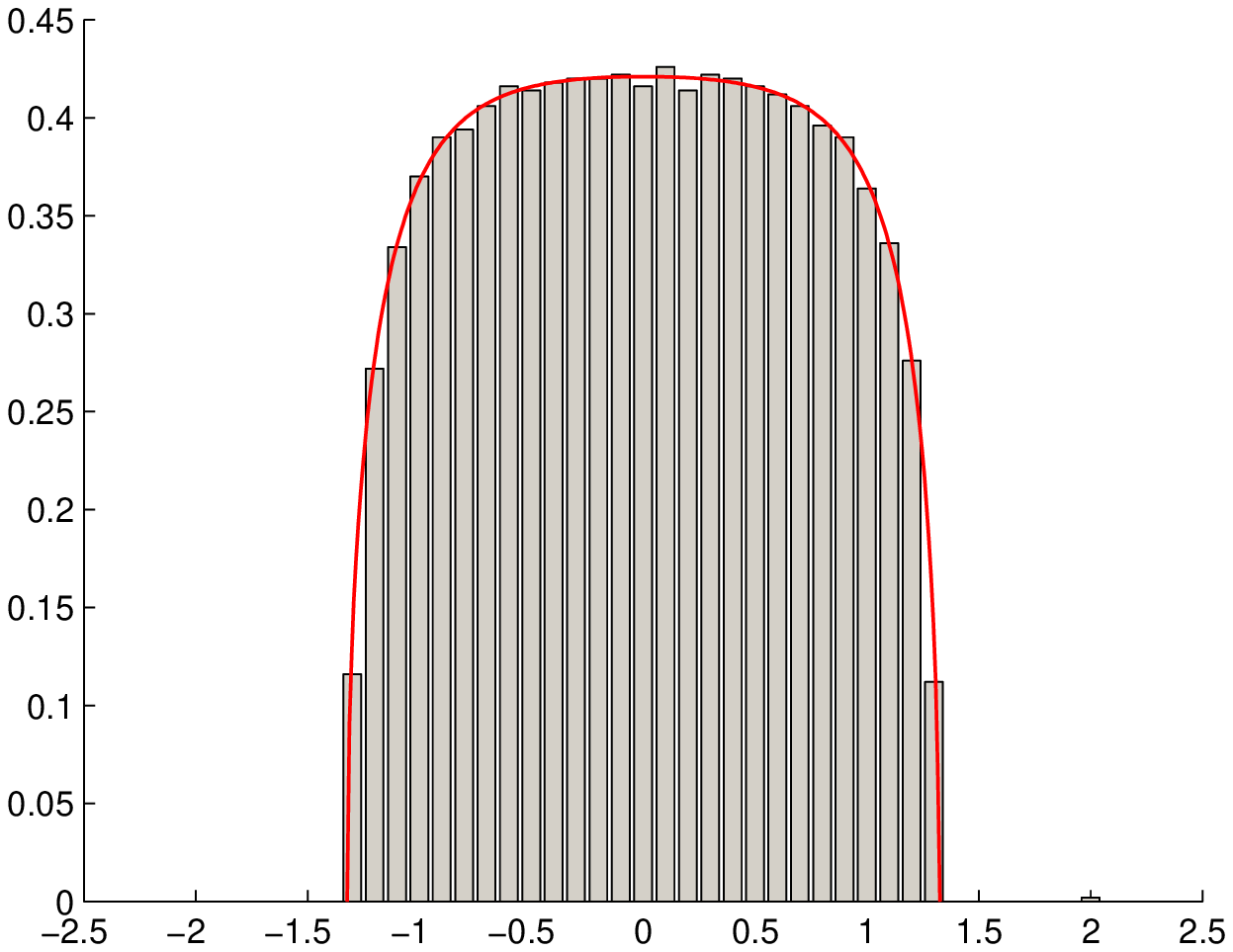}
~~~~~
\includegraphics[width=8cm,height=7cm]{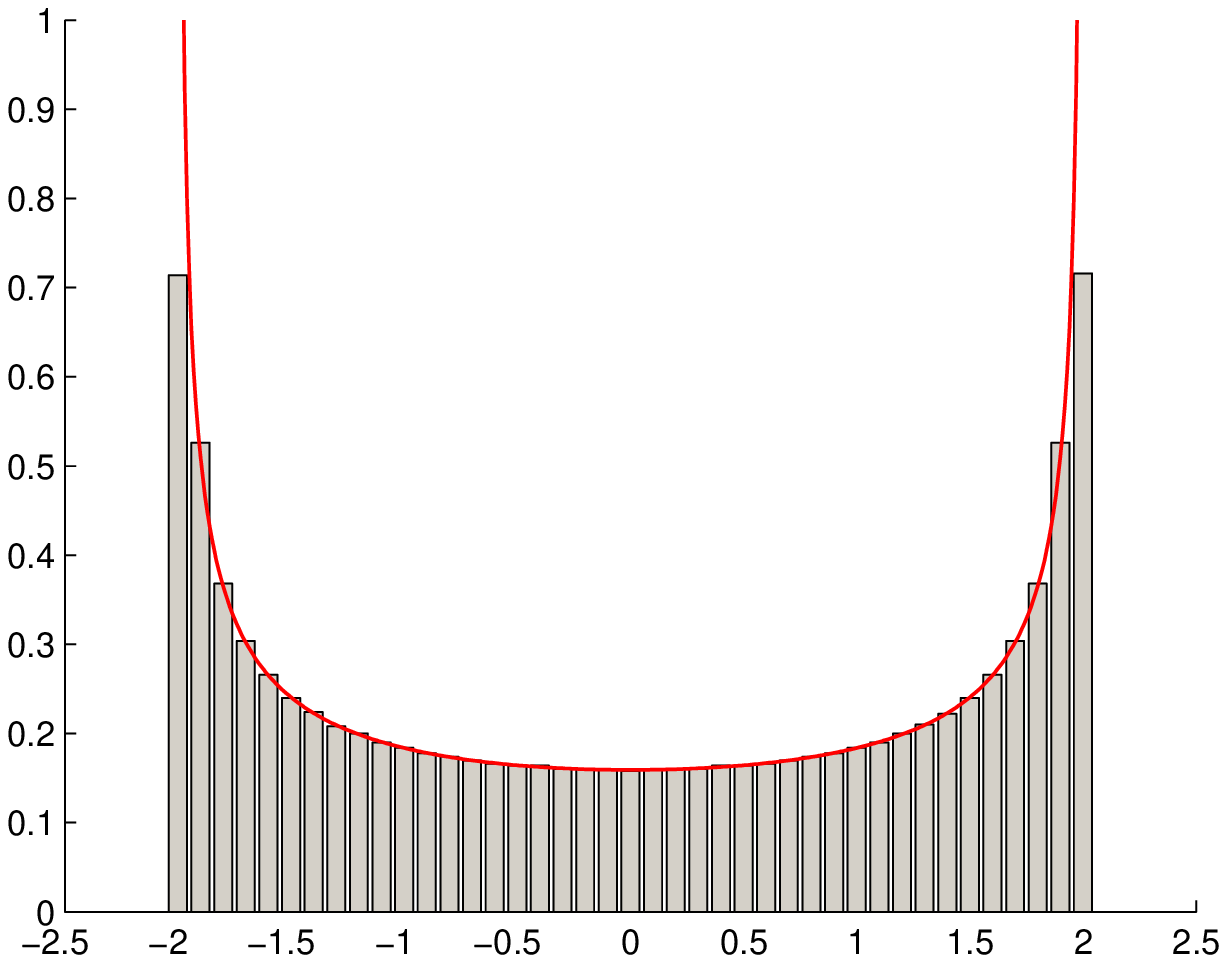}
    \caption{ \label{figJacobi1} \it Density of the limiting distribution
of the random eigenvalues corresponding to the
Jacobi ensemble and a histogram based on $n=5000$ eigenvalues from the Jacobi ensemble (\ref{1.1}). Left panel:  $a_n=b_n=3n,\ \beta_n=2$; right panel:
 $a_n=b_n=\sqrt{n},\ \beta_n=2n$. }
\end{figure}

{\bf Example 3.3.} Assume that
\begin{eqnarray*}
\lim_{n \rightarrow \infty} \frac{\ta}{n}
= \alpha_0 \qquad \mbox{and} \qquad
\lim_{n \rightarrow \infty} \frac{\tb}{n}
= \beta_0
\end{eqnarray*}
for some constants $\alpha_0, \beta_0 \geq 0$. If additionally $(a_n + b_n) / \log n \to \infty$,  it is easy to see that the
assumptions of Theorem 3.2 are satisfied with  $\delta_n = \frac {1}{2}, \ \varepsilon_n=\frac {1}{2} \ (n \in \mathbb{N})$.
Consequently, it follows that the empirical
spectral distribution of the Jacobi ensemble (\ref{1.1}) converges almost surely to a distribution function with density
$$
f_{\alpha_0, \beta_0} (x) = \frac {2 + \alpha_0 + \beta_0}{2 \pi} \ \frac {\sqrt{(2 r_2 - x) (x- 2r_1)}}{4 - x^2} \ I_{(2r_1, 2r_2) } (x) \: ,
$$
where
\begin{align*}
r_1 &:= \frac{\beta_0^2 -\alpha_0^2
-4\sqrt{(\alpha_0+1)(\beta_0+1)(\alpha_0 +\beta_0 +1)}}{(2+\alpha_0
+\beta_0)^2} \ , \\
r_2 &:= \frac{\beta_0^2 -\alpha_0^2
+4\sqrt{(\alpha_0+1)(\beta_0+1)(\alpha_0 +\beta_0 +1)}}{(2+\alpha_0
+\beta_0)^2} \ .
\end{align*}
For example, if $a_n = b_n = 3 n$ and $\beta_n = 2$ we have
\be \label{bild1}
\lim_{n \to \infty} \frac {\tilde a_n}{n} = \lim_{n \to \infty} \frac {\tilde b_n}{n} = 3 \: ,
\ee
and the limiting spectral distribution has the density
$$
f_{3,3} (x) = \frac {4}{\pi} \frac {\sqrt{7/4 - x^2}}{4-x^2} \ I_{[- \frac {\sqrt{7}}{2}, \frac {\sqrt{7}}{2}]} (x) \: ,
$$
which is depicted in the left part of Figure \ref{figJacobi1}.
 Similarly, if $a_n = b_n = \sqrt{n}$ and $\beta_n = 2n$ we have
$$
\lim_{n \to \infty} \frac {\tilde a_n}{n} = \lim_{n \to \infty} \frac {\tilde b_n}{n} = 0 \: ,
$$
and the limiting distribution is given by the arc-sine law on the interval $[-2,2]$ with density
$$
f_{0,0} (x) = \frac {1}{\pi} \frac {1}{\sqrt{4-x^2}} \ I_{(-2,2)} (x) \ ,
$$
displayed in the right part of Figure \ref{figJacobi1}.

\bigskip

{\bf Example 3.4.} If the parameters $\tilde a_n$ and $\tilde b_n$ converge to infinity such that
\begin{align*}
\lim_{n \rightarrow \infty} \frac{\ta}{n} =
 \infty ,
\qquad \lim_{n \rightarrow \infty} \frac{\tb}{n} =
\infty ,
\qquad \lim_{n \rightarrow \infty} \frac{\ta}{\tb} =
\gamma >0 \: ,
\end{align*}
and additionally
\be \label{3.5}
 \beta _n \sim C \qquad \ta = O \ (n^{1+\nu}) \: ,
\ee
for some $C \in \mathbb{R}^{+},\ \nu \in (0,1)$ it follows that the scaled empirical distribution of the eigenvalues of the Jacobi ensemble (\ref{1.1})
with parameters $a_n, b_n$ and $\beta_n$ converges almost
surely to Wigner's semi circle law, that is
\begin{align*}
\lim_{n \rightarrow \infty} F_n^J \left( 2\sqrt{\frac{n}{\ta-1}} \xi -
2\frac{\ta-\tb}{\ta+\tb -2} \right) \stackrel{a.s.}{=}
                    \frac{2}{\pi \sigma^2}
\int_{-\sigma}^{\xi} \sqrt{\sigma^2-x^2} dx \qquad \mid \xi \mid \leq
\sigma \ ,
\end{align*}
where $\sigma = 4 \gamma/ (1+ \gamma )^{3/2}$. Note that assumption (\ref{3.5}) guarantees that the second condition (\ref{3.3}) in Theorem 3.1 is satisfied, i.e.
$$
\lim_{n \rightarrow \infty} \sqrt{\frac{n}{\ta-1}}^4 \frac{a_n
+b_n}{\log n} = \infty .
$$
This situation is also of particular interest in the case,
where the inverse temperature $\beta_n$ converges to 0. For example, if $a_n = b_n = n-1$ and $\beta_n = 2n^{-1/4}$ we have
$\tilde a_n = \tilde b_n =  n^{5/4}$ and the condition (\ref{3.3}) is satisfied. The limiting distribution of the scaled eigenvalues
$$ \frac{1}{2} \sqrt{\frac{n^{5/4}-1}{n}} \lambda^{(n)}_1 < \dots < \frac{1}{2} \sqrt{\frac{n^{5/4}-1}{n}} \ \lambda^{(n)}_n $$
is then given by
$$
f(x) = \frac {1}{\pi} \sqrt{2 - x^2} \ I_{[- \sqrt{2}, \sqrt{2}]} (x) .
$$

\bigskip

{\bf Example 3.5.} We now consider the case, where the sequences $\tilde a_n$ and $\tilde b_n$ converge to infinity with different rates.
If
\begin{align*}
\lim_{n \rightarrow \infty} \frac{\ta}{n} =
\infty \ ,\qquad
\lim_{n \rightarrow \infty} \frac{\tb}{n} =
\beta_0 \geq 0 \: ,
\end{align*}
and the sequences $(\beta _n)_{n \in \mathbb{N}},\ (a_n)_{n \in \mathbb{N}}$ satisfy
$$
\beta _n \sim C , \qquad \ta = O\ (n^{1+\nu})
$$
for constants $C \in \mathbb{R}^{+}$ and $\nu  \in (0,  1/3),$ then it is easy to see that the
conditions of Theorem 3.2 are satisfied with
$$
\varepsilon_n = 0 \: ; \quad
\delta_n = \frac {n}{\ta - 1},
$$
and it follows by a straightforward calculation that
$$
\lim_{n \rightarrow \infty} F_n^J\left( \frac{2n}{\ta-1} \xi -2 \right) \stackrel{a.s.}{=}
\frac{1}{4\pi} \int_{s_1}^{\xi } \frac{\sqrt{(s_2-x)(x-s_1)}}{x} dx
            \qquad s_1
\leq \xi \leq s_2 \: ,
$$
where
\begin{align*}
s_1 &:= 2(2+\beta_0) - 4\sqrt{1+\beta_0 } \ , \\
s_2 &:= 2(2+\beta_0) + 4\sqrt{1+\beta_0 } \ .
\end{align*}
Similarly, if
\begin{align*}
\lim_{n \rightarrow \infty} \frac{\ta}{n} =
\infty ,
\qquad
\lim_{n \rightarrow \infty} \frac{\tilde b_n}{n}= \infty , \qquad \lim_{n \rightarrow \infty} \frac{\ta}{\tb} =
\infty \: ,
\end{align*}
and for some $\nu \in (0,1)$ there exist constants $C, C_1, C_2 > 0$ and a $\mu > 3/2 \nu -1/2$ such that
\be \label{3.6} \beta _n \sim C, \quad \ta \leq C_1n^{1+\nu}, \qquad \tb \geq C_2n^{1+\mu} \: , \ee
then we obtain for any $-2 \leq \xi \leq 6$
$$
\lim_{n \rightarrow \infty} F_n^J\left( 2 \frac{\sqrt{n(\tb-1)}}{\ta-1} \xi - 2
\frac{\ta + 2\sqrt{n(\tb -1)}-\tb}{2n+\ta+\tb-2} \right)
                        \stackrel{a.s.}{=}
\frac{1}{8\pi} \int_{-2}^{\xi } \sqrt{(6-x)(x+2)} dx .
$$

\section{An application to the multivariate $F$ distribution}
\def\theequation{4.\arabic{equation}}
\setcounter{equation}{0}

Let $X_{ij} \ \ (i = 1, \dots, n; \ j=1, \dots, n_1;\ n_1 \geq n)$ and $Y_{ij} \ \ (i = 1, \dots, n; \ j=1, \dots, n_2; \  n_2 \geq n)$ denote independent
standard normal distributed random variables and consider the random matrices \be \label{4.1} X_n = \left ( X_{ij} \right )^{j=1, \dots,
n_1}_{i=1,\dots,n} \in \mathbb{R}^{n\times n_1} \: ; \qquad Y_n = \left ( Y_{ij} \right )^{j=1, \dots, n_2}_{i=1,\dots,n}\in
\mathbb{R}^{n\times n_2} \: . \ee The matrix \be \label{4.2} F_n := \left( \frac{1}{n_1}X_nX_n^T \right) \left( \frac{1}{n_2}Y_nY_n^T
\right)^{-1} \in \mathbb{R}^{n \times n} \ee is called multivariate $F$-matrix and plays a prominent role in the multivariate analysis of
variance [see e.g.\ Muirhead (1982)]. Silverstein (1985b) showed that under the conditions
$$
\lim_{n \to \infty} \frac {n}{n_1} = y > 0 \qquad  \lim_{n \to \infty} \frac {n}{n_2} = y^\prime \in (0,1)
$$
the empirical distribution of the eigenvalues of a multivariate $F$-matrix converges in probability to a non-random distribution function, say $F_{y,y^\prime}$,
with density \be \label{4.3} f_{y,y'}(x) = \frac{1-y'}{2\pi x(xy'+y)} \sqrt{(x-s_1)(s_2-x)} I_{(s_1,s_2)}(x)  \: , \ee where the bounds of the
support are given by
\begin{eqnarray*}
s_1 &=& \Bigl( \frac{1-\sqrt{1-(1-y)(1-y')}}{1-y'} \Bigr) ^2 \ ,\\
s_2 &=& \Bigl( \frac{1+\sqrt{1-(1-y)(1-y')}}{1-y'} \Bigr) ^2 \ .
\end{eqnarray*}
Moreover, if $y > 1$ the limiting distribution has mass $1 - 1/y$ at the point 0. A corresponding result for the expectation of spectral
distribution of the matrix $F_n$ can be found in Collins (2005). In the following discussion we will extend these  results in two different directions
using the methodology developed in Section 2. On the one hand, we prove that in the case $y \in (0,1]$ these results are also correct, if almost sure convergence is
considered. On the other hand, we extend these results to the case where $n, n_1, n_2$ are not necessarily of the same order.

\bigskip

{\bf Theorem 4.1.} {\it Consider the multivariate $F$-matrix defined by (\ref{4.2}). If  \be \label{4.5} \lim_{n \to \infty} \frac {n}{n_1}
= y \in (0, 1] \: , \qquad \lim_{n \to \infty} \frac {n}{n_2} = y^\prime \in (0, 1), \ee then the empirical distribution function
of the eigenvalues of the matrix $F_n$ converges a.s.\ to to a distribution function with density (\ref{4.3}). }

\bigskip

{\bf Proof.} Consider the matrix \be \label{4.6} A_n := 2 \left( Y_nY_n^T-X_nX_n^T \right) \left( Y_nY_n^T+X_nX_n^T \right) ^{-1}  \: , \ee then it
follows from Muirhead (1982) that the joint density of the eigenvalues of $A_n$ is given by the Jacobi ensemble (\ref{1.1}) with $\beta = 1$
and
\begin{eqnarray*}
a_n &=& \frac {1}{2} (n_1 - n - 1)\ , \\
b_n &=& \frac {1}{2} (n_2 - n - 1)\: .
\end{eqnarray*}
If $\lambda^F$ denotes an eigenvalue of the matrix $F_n$ we obtain with some appropriate constant $C \in \mathbb{R}$ the identity
\begin{align*}
& \det \left(  2 \frac{\frac{n_2}{n_1} - \lambda ^F}{\frac{n_2}{n_1} + \lambda ^F} I_n -
                                    2\left( Y_nY_n^T-X_nX_n^T \right) \left( Y_nY_n^T+X_nX_n^T \right) ^{-1} \right) \\
    =& \det \left( \lambda^F I_n - \frac{n_2}{n_1} \left( X_nX_n^T \right) \left( Y_nY_n^T \right)^{-1} \right) \cdot  \det \left( Y_nY_n^T \right)
    \cdot C ~=~0 \: ,
\end{align*}
which shows that \be \label{4.7} \lambda^J = 2 \ \frac {\frac {n_2}{n_1}- \lambda^F}{\frac {n_2}{n_1}+  \lambda^F} \ee is an eigenvalue of the
matrix $A_n$. Consequently,  the empirical distribution function $F^F_n$ of the eigenvalues  of the matrix $F_n$ satisfies the relation \be
\label{4.8} F_n^F(\xi ) \stackrel{a.s.}{=} 1-F_n^J \left( 2 \frac{\frac{n_2}{n_1} - \xi}{\frac{n_2}{n_1} + \xi} \right)  \: , \qquad \forall \ \xi \geq 0   \: , \ee where
$F^J_n$ denotes the empirical distribution function corresponding to the Jacobi ensemble (\ref{1.1}) with  parameter $\beta = 1$ (note that only the case $\xi \geq
0$ is of interest here). We now  use Theorem 3.2 to derive the limiting spectral distribution. For this
purpose we identify
$$
\ta-1= \frac{2a_n+2}{\beta} -1 =n_1-n\ , \qquad \tb -1 = \frac{2a_n+2}{\beta} -1 =n_2-n \ ,
$$
and obtain the limits
\begin{eqnarray*}
\lim_{n \rightarrow \infty} \frac{\ta}{n} &=& \lim_{n \rightarrow \infty} \frac{n_1-n+1}{n} = y^{-1}-1 \geq 0 \ ,\\
\lim_{n \rightarrow \infty} \frac{\tb}{n} &=& \lim_{n \rightarrow \infty} \frac{n_2-n+1}{n} = y'^{-1}-1 >0 \ .
\end{eqnarray*}
Moreover, the assumption
$$
\lim_{n \rightarrow \infty} \frac{a_n+b_n}{\log n} = \infty
$$ is obviously satisfied. From  Example 3.3  it therefore follows that the empirical distribution function $F^J_n$ of the eigenvalues of the
matrix $A_n$ converges a.s.\ to a distribution function $F$ with density
$$
f(x) = \frac{y^{-1} +y'^{-1} }{2\pi} \frac{\sqrt{(2r_2-x)(x-2r_1)}}{4-x^2} I_{(2r_1,2r_2)}(x)  \: ,
$$
where $r_1$ and $r_2$ are given by
\begin{eqnarray*}
r_1 &:= &
\frac{(y-yy')^2 -(y'-yy')^2-4\sqrt{y^2y'^2(y+y' -yy')}}{(y+y')^2} \ ,\\
r_2 &:= &\frac{(y-yy')^2 -(y'-yy')^2+4\sqrt{y^2y'^2(y+y' -yy')}}{(y+y')^2}
\end{eqnarray*}
Observing the relation (\ref{4.8})
we obtain
$$
F^F_n (\xi)  \stackrel{a.s.}{\longrightarrow} 1 - F \left ( 2 \frac {y - y^\prime \xi}{y + y^\prime \xi} \right )  \: ,
$$
and the assertion of the theorem now follows by a straightforward but tedious calculation of the density of the limiting
distribution. \hfill $\Box$

\bigskip

While the preceding theorem essentially provides an alternative proof of the results of Silverstein (1985b), the following
three theorems extend Silverstein's findings to the case where $y, y^\prime = 0$.

\bigskip

{\bf Theorem 4.2.} {\it Consider the multivariate $F$-matrix defined in (\ref{4.2}) and denote by $\lambda^F_1 \leq  \dots \leq \lambda^F_n$ the
corresponding eigenvalues. If
$$
\lim_{n \rightarrow \infty} \frac{n}{n_1} = 0 \ , \qquad \lim_{n \rightarrow \infty} \frac{n}{n_2} = 0\ ,
                                        \qquad \lim_{n \rightarrow \infty} \frac{n_1}{n_2} = \gamma >0
$$
and
$$
n_1 = O(n^{1+\nu})
$$
with $\nu \in (0,1)$, then  the empirical distribution function of the transformed eigenvalues
$$
\mu_i = 2 \sqrt{\frac {n_1}{n}-1}  \ \left \{ \frac {n_2 - n}{n_1 + n_2 - 2n} - \frac {n_2}{n_1 \lambda^F_{i} + n_2}  \right \}  \qquad i=1, \dots, n
$$
converges a.s.\ to a distribution function with density
$$
f_\gamma (x) ~=~
\frac{2}{\pi \sigma^2}
 \sqrt{ \sigma^2 - x^2 }        \ I\left \{ -\sigma < x < \sigma \right \} \ ,
$$
where $\sigma = 4 \gamma / (1 + \gamma )^{3/2}$.}

\bigskip

{\bf Proof.} Recall the definition of the matrix $A_n$ in (\ref{4.6}), which corresponds to the Jacobi ensemble (\ref{1.1}) with $\beta = 1,\
a_n = \frac {1}{2} (n_1 - n - 1),\ b_n = \frac {1}{2} (n_2 - n - 1)$. Using the notation (\ref{2.25}) we obtain  $\tilde a_n - 1 = n_1 - n,\
\tilde b_n - 1 = n_2 - n$. By the assumption of the theorem we have
\begin{align*}
\lim_{n \rightarrow \infty} \frac{\tilde a_n}{n} = \infty \ ,  \qquad \lim_{n \rightarrow \infty} \frac{\tilde b_n}{n} = \infty \ , \qquad
\lim_{n \rightarrow \infty} \frac{\tilde a_n}{\tilde b_n} = \gamma
\end{align*}
and $\ta = O (n^{1+\nu})$. Therefore it follows from  Example 3.4 that
\begin{align*}
& \lim_{n \rightarrow \infty} F_n^J \left( 2\sqrt{\frac{n}{n_1-n}} \xi - 2\frac{n_1-n_2}{n_1+n_2-2n} \right) \\
= &\frac{2}{\pi \sigma^2} \int_{\infty}^{\xi} \sqrt{\sigma^2 -x^2 }
                    I\left \{ -\sigma < x < \sigma \right \} dx
\end{align*}
a.s., where $F^J_n$ denotes the empirical distribution of the eigenvalues of the matrix $A_n$. The identity (\ref{4.8}) implies for $\xi > - 2$
$$
F_n^J(\xi ) \stackrel{a.s.}{=} 1-F_n^F \left( \frac{n_2}{n_1} \frac{2 - \xi}{2 + \xi} \right) \qquad
$$
and therefore it follows

\begin{align*}
\lim_{n \rightarrow \infty} F_n^F \left( \frac{n_2}{n_1} \left( \frac{2}{2\frac{n_2-n}{n_1+n_2 - 2n} - \sqrt{\frac{n}{n_1-n}} \xi } -1 \right) \right)
= \frac{2}{\pi \sigma^2} \int_{\infty}^{\xi} \sqrt{\sigma^2 -x^2 }
                    I\left \{ -\sigma < x < \sigma \right \} dx
\end{align*}
a.s., which proves the assertion of the theorem. \hfill $\Box$

\bigskip

{\bf Theorem 4.3.} {\it Consider the multivariate $F$-matrix and denote by $\lambda^F_1 \leq \dots \leq \lambda^F_n$ the
corresponding eigenvalues.
If
$$
\lim_{n \rightarrow \infty} \frac{n}{n_1} = 0\ , \qquad \lim_{n \rightarrow \infty} \frac{n}{n_2} = y' \in (0,1]
$$
and
$$
n_1 = O(n ^{1+\nu}) \: ,
$$
with $\nu \in (0,1/3)$, then the empirical distribution function of the scaled eigenvalues
$$
\mu_i = \frac {n}{2 (n_1 - n)} \left ( \lambda^F_i \frac {n_1}{n_2} + 1 \right ) \qquad i=1, \dots, n
$$
converges a.s.\  to a distribution function $F$ with density
$$
f_{y'} (x) = \frac{1}{4\pi} \frac{\sqrt{(xs_2-1)(1-xs_1)}}{x^2} I_{(s_2^{-1},s_1^{-1})} (x)  \: ,
$$
where the bounds of the support of the density are given by}
\begin{align*}
s_1 &:= 2(y'^{-1}+1) - 4\sqrt{y'^{-1}} \ ,\\
s_2 &:= 2(y'^{-1}+1) + 4\sqrt{y'^{-1}} \ .
\end{align*}

\bigskip

The proof is analogous to the proof of Theorem 4.2, using the first result of Example 3.5.
Similarly, the following theorem can be proven using the last statement in Example 3.5.

\bigskip

{\bf Theorem 4.4.} {\it Consider the multivariate $F$-matrix and denote by $\lambda^F_1 \leq \dots \leq \lambda^F_n$ the
corresponding eigenvalues.
If
$$
\lim_{n \rightarrow \infty} \frac{n}{n_1} = 0\ , \qquad \lim_{n \rightarrow \infty} \frac{n}{n_2} = 0,\ \qquad
                                                                                                                \lim_{n \rightarrow \infty} \frac{n_1}{n_2} = \infty
$$
and for some $\nu \in (0,1)$ there exist constants $C_1, C_2 > 0$ and a $\mu > 3/2 \nu -1/2$ such that
\begin{equation*} n_1 \leq C_1n^{1+\nu}, \qquad n_2 \geq C_2n^{1+\mu},\ \end{equation*}
then the empirical distribution function of the scaled eigenvalues
$$
\mu_i = 2\frac {n_1-n}{n_1+n_2} \frac{\frac{n_1}{n_2} (n_2 - \sqrt{n(n_2-n)}) \lambda^F_i - (n_1+\sqrt{n(n_2-n)})}
                                                                            {\sqrt{n(n_2-n)} (1+\frac{n_1}{n_2} \lambda^F_i)}  \qquad i=1, \dots, n
$$
converges a.s.\  to a distribution function $F$ with density
$$
f (x) = \frac{1}{8\pi} \sqrt{(6+x)(2-x)} I_{[-6,2]}(x) \ .
$$}

\section{Appendix: auxiliary results}
\def\theequation{5.\arabic{equation}}
\setcounter{equation}{0}

{\bf Lemma A.1.} {\it Let $Z$ denote a Beta-distributed random variable on the interval  $[0,1]$ with density \be \label{a.1} \frac {\Gamma
(p+q)}{\Gamma(p) \Gamma(q)} \ x^{p-1} (1-x)^{q-1} \ I_{(0,1)} (x) \quad \quad (p,q>0) \: , \ee then for any $\delta > 0$
$$
P ( \ \mid Z - E [Z] \mid \ > \delta ) \leq 4 e^{c(p+q)}  \: ,
$$
where the constant $c$ is defined by
$$
 c=  \log \left( 1+\frac{\delta }{3+2\delta } \right) - \frac{\delta }{3+2\delta } \ .
$$
}

\bigskip

{\bf Proof.} If $X \sim \Gamma (p, \ p+q), \ Y \sim \Gamma (q, \ p+q)$ denote independent Gamma-distributed random variables, it is well known
that  the ratio $Z=X/(X+Y)$ has a Beta-distribution with density (\ref{a.1}). Because $E[Z] = E[X] = p/(p+q)$ it follows that \be \label{a.2}
P(|Z - \erw{Z} | > \delta ) = P\Bigl( \Bigl| \frac{X}{X + Y} - \erw{X} \Bigr| > \delta \Bigr) \ee Define $\delta^\prime = \delta / (3 + 2
\delta)$ and assume that
$$
\mid X - E[X] \mid \ \leq \delta^\prime \: , \quad \mid Y - E [Y] \mid \ \leq \delta^\prime \: ,
$$
then it is easy to see that
$
\mid X + Y - 1 \mid \ = \ \mid X + Y - E [X+Y] \mid \ \leq 2 \delta^\prime \ ;
$
and
$$
\Bigl| \frac{X}{X + Y} - \erw{X} \Bigr| \leq  \frac{1}{1-2\delta^\prime} ( | X - \erw{X} | + \erw{X} | X + Y - \erw{X+ Y} | ) \leq
\frac{1}{1-2\delta^\prime}  3\delta^\prime = \delta \: .
$$
This implies for the probability in (\ref{a.2})
\begin{eqnarray} \label{a0}
  P\left( \left| \frac{X}{X+ Y} - \erw{X} \right| > \delta \right) &\leq &
P(X > \erw{X} + \delta^\prime ) + P(X< \erw{X} - \delta^\prime )  \\ \nonumber
&& ~~~~~~~~~ + P(Y> \erw{Y} + \delta^\prime ) + P(Y < \erw{Y} +
\delta^\prime ) \ .
\end{eqnarray}
Using similar arguments as in Dette and Imhof (2007) we obtain the estimates
\begin{eqnarray*}
P ( U > E [U] + \delta^\prime) &\leq& \exp \{ (p+q) (\log (1+\delta^\prime) - \delta^\prime) \} \ , \\
P ( U < E [U] - \delta^\prime) &\leq& \exp \{ (p+q) (\log (1-\delta^\prime) + \delta^\prime) \} \ ,
\end{eqnarray*}
where $U$ is either $X$ or $Y$. The assertion of Lemma A.1 now follows from (\ref{a0}) and the definition of  $\delta^\prime$ observing that $
\log (1+ \delta^\prime) - \delta^\prime > \log (1 - \delta^\prime) + \delta^\prime $. \hfill $\Box$

\bigskip
\medskip

{\bf Acknowledgements.} The authors are grateful to Martina Stein  who typed parts of this paper with
considerable technical expertise.
The authors would also like to thank two anonymous referees for their constructive
 comments on an earlier version of
this paper.
The work of the authors was supported by the Sonderforschungsbereich Tr/12
 (project C2, Fluctuations and universality of invariant random matrix ensembles) and
in part by a NIH grant award
IR01GM072876:01A1.

\bigskip
\medskip
\medskip

{\large References}

M. Abramovich, I. Stegun (1965). Handbook of Mathematical Functions. Dover Publications Inc., N.Y.

\smallskip

Z. D. Bai (1999). Methodologies in spectral analysis of large dimensional
random matrices, a review. Statist.\ Sinica 9, 611-677.

\smallskip

T. S. Chihara (1978). An Introduction to Orthogonal Polynomials. Gordon and Breach, New York.

\smallskip

B. Collins (2005). Product of random projections, Jacobi ensembles and universality problems arising from free probability.
    Probability Theory and Related Fields 133, 315-344

\smallskip

{H. Dette and L. A. Imhof (2007). Uniform approximation of eigenvalues in Laguerre and Hermite $\beta$-ensembles
by roots of orthogonal polynomials. Transactions of the American Mathematical Society 359, 4999-5018.}

\smallskip

H. Dette and W.J. Studden (1995). Some new asymptotic properties
 for the zeros of the Jacobi, Laguerre and Hermite polynomials.
 Constructive Approximation, Vol. 11, 227-238.

 \smallskip

I. Dumitriu and A. Edelman (2002).  Matrix models for beta-ensembles, J. Math. Phys. 43, 5830-5847.

\smallskip

F. J. Dyson (1962). The threefold way. Algebraic structure of symmetry groups and ensembles
in quantum mechanics. J. Math.\ Phys. 3, 1199-1215.

 \smallskip

 Edelman and Sutton (2006). The beta-Jacobi matrix model, the CS decomposition, and generalized singular value problems.
 Found. Comput. Math. In press. http://dx.doi.org/10.1007/s10208-006-0215-9.

\smallskip

A. Elbert, A. Laforgia, L. G. Rodono (1994). On the zeros of Jacobi polynomials. Acta Math. Hungar. 64, 351-359.

\smallskip

W. Gawronski, B. Shawyer (1991). Strong asymptotics and the limit distribution of the zeros of Jacobi polynomials
$P_n^{(a_n + \alpha, b_n +\beta)} (x)$. In: Progress in Approximation Theory (P. Nevai, A. Pinkus eds.). New York: Academic Press 379-404.

 \smallskip

{R. Killip  and I. Nenciu I (2004). Matrix models for circular ensembles. Int. Math. Res. Not. 50, 2665-2701.}

%

\smallskip
A. B. J. Kuijlaars and W. van Assche (1999). The asymptotic zero distribution of orthogonal polynomials with varying recurrence coefficients,
J. Approx. Theory 99 (1999), 167-197.

\smallskip

R. A. Lippert (2003). A matrix model for the beta-Jacobi ensemble. J. Math. Phys. 44 4807-4816.

\smallskip

 R. J. Muirhead (1982). Aspects of Multivariate Statistical Theory. Wiley, New York.

\smallskip

J. W. Silverstein (1985a) The smallest eigenvalue of a large
dimensional Wishart matrix,
 Ann.\ Probab.\ 13, 1364-1368.

\smallskip

J. Silverstein (1985b). The limiting eigenvalue distribution of a multivariate $F$-matrix. SIAM J. Math. Anal. 16, 641-646.

\smallskip

G. Szeg\"o (1975). Orthogonal Polynomials. 4th ed. Amer.\ Math.\ Soc., Providence, RI.

\end{document}